\documentclass[11pt,twoside]{article}
\usepackage{amsmath}
\usepackage{amsthm}


\usepackage{mathtools}
\usepackage{amssymb}
\usepackage{graphicx}
\usepackage{enumerate}
\usepackage[authoryear,sort&compress,round,comma]{natbib}

\usepackage{mathptmx}

\newtheorem{theorem}{Theorem}[section]
\newtheorem{lemma}{Lemma}[section]

\numberwithin{equation}{section}

\usepackage{graphicx,type1cm,eso-pic,color}
\makeatletter
\AddToShipoutPicture{
	\setlength{\@tempdimb}{.5\paperwidth}
	\setlength{\@tempdimc}{.5\paperheight}
	\setlength{\unitlength}{1pt}
	\put(\strip@pt\@tempdimb,\strip@pt\@tempdimc)
}
\makeatother

\begin{document}
	\setcounter{page}{1}

	\thispagestyle{empty}
	\markboth{}{}

	\pagestyle{myheadings}
	\markboth{}{ }
	
	\pagestyle{myheadings}
	\markboth{N.Gupta and S.K.Chaudhary}{N.Gupta and S.K.Chaudhary}
	\date{}
	
	
	\noindent 
	
	\vspace{.1in}
	
	{\baselineskip 20truept
		
		\begin{center}
			{\Large {\bf Some characterizations of continuous symmetric distributions based on extropy of record values }} \footnote{\noindent
			{\bf *} E-mail: nitin.gupta@maths.iitkgp.ac.in\\
			{\bf ** }  corresponding author E-mail: skchaudhary1994@kgpian.iitkgp.ac.in}\\
		
	\end{center}
	
	\vspace{.1in}
	
	\begin{center}
		{\large {\bf Nitin Gupta* and  Santosh Kumar Chaudhary**}}\\
		{\large {\it Department of Mathematics, Indian Institute of Technology Kharagpur, West Bengal 721302, India. }}
		\\
	\end{center}
}
	\vspace{.1in}
	\baselineskip 12truept

	\begin{center}
		{\bf \large Abstract}\\
	\end{center}
     Using different extropies of $k$ record values various characterizations are provided for continuous symmetric distributions. The results are in addition to the results of Ahmadi, J. (Statistical Papers, 2021, 62:2603-2626). These include cumulative residual (past) extropy, generalised cumulative residual (past) extropy, also some common Kerridge inaccuracy measures. Using inaccuracy extropy measures, it is demonstrated that continuous symmetric distributions are characterised by an equality of information in upper and lower $k$-records.

	\section{Introduction}
	 Let $X$ be a continuous Random Variable with probability density function (pdf) $f_X(x)$, cumulative distribution function (cdf) $F_X (x)$ and support $S_X$. A probability distribution is said to be symmetric around $k$ if and only if there exists a finite number $k$ such that $f_X(k+x)=f_X(k-x) $ for all $x\in S_X$  that is, $F_X(k-x)+F_X(k+x)=1$ for all $x\in S_X$. Normal distribution, Laplace distribution, Cauchy distribution, Uniform distribution, Logistic distribution, and Student's t-distribution are some well-known symmetric distributions. See Johnsons et al. (1995) for a variety of examples of characterization and application related to symmetric distribution. In probability and statistics, symmetry is a key structural presumption that is relevant to a wide range of issues. In finance, the Sharpe-Lintner capital asset pricing model and the Black–Scholes option pricing model are widely used models that need symmetry as an assumption. Various features of the symmetry of probability distributions have been extensively studied in the literature. Properties of symmetric distribution have been studied by several authors. Traders use symmetric distribution in the financial sector to predict the value of a stock, currency or commodity in a time frame. For example, Based on a symmetric distribution of price over time, the mean reversion hypothesis in finance predicts that asset prices would eventually return to their long-term mean or average values, see Ahmed et al. (2018).	
	 																																																																																																																																																																 \indent Based on the characteristics of information measures of record data obtained from a sequence of independently and identically distributed (iid) continuous random variables, Ahmadi (2020) gave some new characterizations for continuous symmetric distributions. A nonparametric estimator for the symmetric distribution function under multi-stage ranked set sampling was proposed by Mahdizadeh and Zamanzade (2020). Usually, Based on the distinct characteristics of a symmetric distribution, criteria for determining whether a distribution is symmetric or not can be developed. Hence, a goodness-of-fit test for symmetry can be constructed using characterization results, see, for example, Dai et al. (2018) and Bozin et al. (2020) and references therein. Characterizations of symmetric distribution are given by many researchers, see Ahmadi and Fashandi (2019), Ahmadi(2019, 2020, 2021), Ushakov (2011), and Ahmadi et al. (2020). One of the possible applications of the characterization could be testing for symmetry. Xiong et al. (2021) proposed a test for symmetry using extropy of upper $k$-record value and lower $k$ record value. Jose and Sather (2022) proposed a test for symmetry using extropy of $n$th upper $k$-record value and $n$th lower $k$ record value. The primary objective of this work is to provide some new, extropies-based characterizations for continuous symmetric distributions.
	
	 Suppose that a random variable $X$  has cdf $F_X(x)$ with continuous  pdf $f_X(x) $. To measure	the uncertainty contained in $X$, the entropy was defined by Shannon (1948) as follows, $$H(X)=-\int_{S_X} f_X(x) \log f_X(x)dx.$$
	 The extropy of the random variable X is defined by Lad et al. (2015) to be
	$$J(X)=-\frac{1}{2} \int_{S_X} f_X^2(x)dx.$$ Lad et al. (2015) described the characteristics of this measure, including the maximum extropy distribution and statistical applications. Also, Useful results can be found in Qiu (2017) and Qiu and Jia (2018) related to extropy and residual extropy properties of order statistics and record values. For more details on extropy see Raqab and Qiu (2019), Noughabi and Jarrahiferiz (2019), Krishnan et al. (2020), Balakrishnan et al. (2020), Bansal and Gupta (2021), Qiu and Raqab (2022), Gupta and Chaudhary (2022) and the references therein.

	\indent The concept of $k$-records was introduced by Dziubdziela and Kopociński (1976); for more details, also see  Ahsanullah (1995) and Arnold et al. (1998).
	 The pdf of the $n$-th upper $k$-record value $U_{n,k}$ and the $n$th lower $k$-record value $L_{n,k}$ respectively are given by (see Arnold et al.(2008) and  Ahsanullah (2004)) 
	 
	 	\begin{align*}
	 	&f_{U_{n(k)}}(x)=\dfrac{k^n}{(n-1)!}[-\log \overline{F}(x)]^{n-1} \overline{F}(x)^{k-1}f_X(x),\label{r1}\\
	 	\text{and} \ \ \ &g_{L_{n(k)}}(x)=\dfrac{k^n}{(n-1)!}[-\log F_X(x)]^{n-1} F_X(x)^{k-1}f_X(x).\\ 
	 \end{align*}
	 
	\noindent The CDF of $U_{n,k}$ and  $L_{n,k}$, respectively, are  
	\begin{align*}
		F_{U_{n,k}}(u)&=  1-\bar{F}^k_X (u)   \sum_{i=0}^{n-1} \frac{(-k \log\bar{F}_X (u) )^i}{i!},\\
	\text{and} \ \ \ 	F_{L_{n,k}}(u)&=  {F}^k_X (u)   \sum_{i=0}^{n-1} \frac{(-k \log F_X (u) )^i}{i!}.
	\end{align*}
	
    \noindent The extropy of random variable $U_{n,k}$ and $L_{n,k}$, respectively, are 
	
		\begin{align}
		&J(U_{n,k}) =-\frac{1}{2} \int_{0}^{\infty} {f}_{U_{n,k}}^2(x)dx, \\ 
		\text{and} \ \ \  &J(L_{n,k}) =-\frac{1}{2} \int_{0}^{\infty} {f}_{L_{n,k}}^2(x)dx.
	    \end{align}
	
	\noindent Xiong et el. (2021) showed that the pdf of $X$ is symmetric about its finite mean $\mu$ if and only if $J(U^X_k)=J(L^X_k)$ for all $k=1,2,3,\dots$ where $U^X_k$ and $L^X_k$ represent $k$th upper record value and $k$th lower record value, see Arnold et al.(1998) for more on upper and lower record value. Jose and Sathar (2022) showed that $X$ is symmetric about its finite mean $\mu$ if and only if  $J(L_{n,k})=J(U_{n,k})$ for all $k=1,2,3,...$. Ahmadi (2021) provided several characterizations of continuous symmetric distribution based on information measures of $k$-record.
			
    In this paper, we present some characteristics of the symmetric distribution. We discuss our results in section 2 using cumulative residual (past) extropy. We provide the results in section 3 using an inaccuracy measure. We used a few examples in section 4 to demonstrate the results of the previous sections. Section 5 concludes with a suggestion for potential further research.

	\section{Results based on cumulative residual (past) Extropy  }
	
	The following lemma due to Fashandi and Ahmadi (2012) will be used to develop different characterization of symmetric distribution.
	
	\begin{lemma} \label{fa2012}(Fashandi and Ahmadi, 2012) Let $X$ be a continuous	random variable with pdf $f_X$ and cdf $F_X$ with support $S_X$. Then, the	identity,  $$f_X(F_X^{-1} (u))=f_X(F_X^{-1} (1-u))$$ for almost all $ u \in (0,\frac{1}{2})$ if and only if  that there exists a constant $k$ such that $F_X(k-x)+F_X(k+x)=1$ for all $x\in S_X.$
	\end{lemma}

  Let $\mathbb{C}$ denote the class of all continuous pdf $f_X$, having cdf $F_X$ such that  $$f_X\left(F_X^{-1}(1-u)\right)\geq(\leq)f_X\left(F_X^{-1}(u)\right)$$ for all $u\in (0,\frac{1}{2})$. It can be observed using Lemma \ref{fa2012} that $F$ is symmertic if and only if $f_X(F_X^{-1} (u))=f_X(F_X^{-1} (1-u))$ for almost all $ u \in (0,\frac{1}{2})$. The class $\mathbb{C}$  is non empty and includes various distributions. One may refer Ahmadi et al. (2020) for these distributions.\\
  	\indent Cumulative residual entropy (CRE) was introduced by Rao et al. (2004) in terms of the survival function of $X$. Based on the cdf of $X$, the cumulative past entropy (CPE) was proposed by Di Crescenzo and Longobardi (2009). 	Cumulative residual extropy (CRJ) was introduced by Jahanshahi et al (2020) and it is defined as
	\begin{align}\label{a}
	\xi J(X) &=-\frac{1}{2} \int_{0}^{\infty} \bar{F_X}^2(x)dx .
	\end{align}
    Cumulative past extropy (CPJ) is defined as
    \begin{align}\label{b}
    \bar{\xi} J(X) &=-\frac{1}{2} \int_{0}^{\infty} {F_X}^2(x)dx. 
    \end{align}
The following theorem gives a relationship of symmetric distribution with CPJ and CRJ.
  \begin{theorem} \label{g}  The following two statements are equivalent for any $F_X\in \mathbb{C}:$
 \begin{enumerate}[(i)] \item random variable $X$ has a symmetric distribution;
\item $\bar{\xi} J(X)=\xi J(X)$. \end{enumerate}  
	\end{theorem} 
\noindent \textbf{Proof}   	From (\ref{a}), the CRJ can be written as
     \begin{align}\label{c}
      \xi J(X) &=-\frac{1}{2} \int_{0}^{1} \frac{u^2}{f_X(F_X^{-1} (1-u))}du, 
     \end{align}
   and from (\ref{b}) the CPJ can be written as
     	\begin{align}\label{d}
      	\bar{\xi} J(X) &=-\frac{1}{2} \int_{0}^{1} \frac{u^2}{f_X(F_X^{-1} (u))}du.
     \end{align}
     When $\bar{\xi} J(X)= \xi J(X)$ holds then from (\ref{c}) and (\ref{d}), we have 
\begin{equation}\label{ddd}
\int_{0}^{1} \eta(u) u^2 du =0,
\end{equation}
     \begin{align}
     \text{where }\ \	 \eta(u)=\frac{1}{f_X(F_X^{-1} (1-u))}-\frac{1}{f_X(F_X^{-1} (u))}.
     \end{align}
   Then as $\eta(1-u)=-\eta(u)$, we can write (\ref{ddd}) as
   \begin{align} \label{e}
   	  \int_{0}^{\frac{1}{2}} \eta(u) (2u-1)du=0 .
   \end{align}
  Note that $2u-1  \leq  0 $ for all $0\leq u \leq \frac{1}{2}$. Since by assumption $F_X\in \mathbb{C}$, hence from (\ref{e}) we have $f_X(F_X^{-1} (u))=f_X(F_X^{-1} (1-u))$, so lemma \ref{fa2012} completes the proof. $\blacksquare$

   

    The cdf of $n$th upper $k$-record value $U_{n,k}(u)$ is 
     \begin{align}\label{unkcdf}
     	F_{U_{n,k}}(u)&=  1-\bar{F}^k_X (u)   \sum_{i=0}^{n-1} \frac{(-k \log\bar{F}_X (u) )^i}{i!},
     \end{align}
 and  the cdf of $n$th lower $k$-record value $L_{n,k}(u)$ is  
 \begin{align}\label{lnkcdf}
 	F_{L_{n,k}}(u)&=  {F}^k_X (u)   \sum_{i=0}^{n-1} \frac{(-k \log F_X (u) )^i}{i!}.
     \end{align}
     We establish the result in theorem  \ref{g} based on CRJ and CPJ of records. The following theorem gives a relationship of symmetric distribution with CPJ and CRJ of $k$-record values.
     \begin{theorem}\label{thm2}  Let $X_1,X_2,\ldots$ be a random sample of continuous random variables from a population $X$ having cdf $F_X$ and pdf $f_X$.
Then following two statements are equivalent for any $F_X \in \mathbb{C}:$
\begin{enumerate}[(i)]   	
    \item   random variable $X$ has a symmetric distribution;
 \item for a fixed $k\geq 1$,  $\bar{\xi} J(L_{n,k})=\xi J(U_{n,k})$ for all $n\geq 1$.
\end{enumerate}
      \end{theorem}
     \noindent \textbf{Proof} \
      From (\ref{a}) and (\ref{unkcdf}), the CRJ of $U_{n,k}$ is 
      \begin{align}\label{crex}
      \xi J(U_{n,k}) &=-\frac{1}{2}\int_{0}^{\infty}\bar{F}_{U_{n,k}}^2(x)dx \nonumber \\ 
      &=-\frac{1}{2} \int_{0}^{1} \left( \sum_{i=0}^{n-1} \frac{(-k \log u)^i}{i!} \right) ^2 \frac{u^{2k}}{f_X(F_X^{-1}(1-u)))}du      
    \end{align}
From (\ref{b}) and (\ref{lnkcdf}), the CPJ of $L_{n,k}$ is 
\begin{align}\label{cpex}
	\bar{\xi} J(L_{n,k}) &=-\frac{1}{2} \int_{0}^{\infty}{F}_{L_{n,k}}^2(x)dx \nonumber \\
	&=-\frac{1}{2} \int_{0}^{1} \left( \sum_{i=0}^{n-1} \frac{(-k \log u)^i}{i!} \right) ^2 \frac{u^{2k}}{f_X(F_X^{-1}(u)))}du.
      \end{align}
      
      When $\bar{\xi} J(L_{n,k})=\xi J(U_{n,k})$ then from (\ref{crex}) and (\ref{cpex}) we get,
      \begin{align}\label{z}
      	 \int_{0}^{1} \phi_n^2(u) \eta(u)du=0,       
      \end{align}
     
     \begin{align*}
      	\text{where} \ \ \  \eta(u)=\frac{1}{f_X(F_X^{-1} (1-u))}-\frac{1}{f_X(F_X^{-1} (u))},\\ \text{and} \ \ 
      	\phi_n(u)= u^k\sum_{i=0}^{n-1} \frac{(-k \log u)^i}{i!}.
      \end{align*}
     Further simplification of (\ref{z}) gives
    \begin{equation}\label{lnx} \int_{0}^{\frac{1}{2}} [\phi_n^2(u)-\phi_n^2(1-u)] \eta(u)du =0. \end{equation}
Note that $\phi_n(u)\geq 0,\ u\in(0,1)$. \label{phinon} Also
\begin{align}
\frac{d\phi_n(u)}{du} \nonumber
&=ku^{k-1}\sum_{i=0}^{n-1} \frac{(-k \log u)^i}{i!}-ku^{k-1}\sum_{i=1}^{n-1} \frac{(-k \log u)^{i-1}}{(i-1)!}\\ \nonumber
&=ku^{k-1}\left(\sum_{i=0}^{n-1} \frac{(-k \log u)^i}{i!}-\sum_{i=0}^{n-2} \frac{(-k \log u)^i}{i!}\right)\\ \nonumber
&=ku^{k-1} \frac{(-k \log u)^{n-1}}{(n-1)!}\\
&\geq 0. \label{phiincresing}
\end{align}    
Since by assumption $F_X\in \mathbb{C}$ and $\phi_n(u)$ is positive and increasing in $u$, hence from (\ref{lnx}) we have $f_X(F_X^{-1} (u))=f_X(F_X^{-1} (1-u))$, so the lemma \ref{fa2012} completes the proof. $\blacksquare$


      Generalized cumulative residual entropy (GCRE) of order m was proposed by Psarrakos and Navarro (2013) which is related to the upper record values of a sequence of iid random variables and with the relevation transform. Generalized cumulative past entropy (GCPE) of order m was introduced by Kayal (2016) as an extended version of cumulative past entropy which is related to the lower records and the reversed relevation transform. \\
   \indent The Generalized cumulative past extropy (GCPJ)  $	G\bar {\xi}_m (J(X))$ and  Generalized Cumulative Residual extropy (GCRJ) $G\xi_m(J(X))$ are, respectively, given by
    \begin{align}
    	G\bar {\xi}_m (J(X))&=-\frac{1}{2}\int_{0}^{\infty}{F}^m_X (x)dx,\label{gxp}\\
    	\text{and }\ \ 
    	G\xi_m(J(X))&=-\frac{1}{2}\int_{0}^{\infty}\bar{F}^m_X (x)dx, \label{gxr}
    	 \end{align}
where $m\geq 1$. The following theorem gives relationship of symmetric distribution with GCPJ and GCRJ.
   \begin{theorem}\label{thm3}   The following two statements are equivalent for any $F_X\in \mathbb{C}$:
 \begin{enumerate}[(i)] \item random variable $X$ has a symmetric distribution;
\item $G\bar {\xi}_m (J(X))=G\xi_m(J(X))$ for a given fixed number $m \geq$ 1. \end{enumerate}  
	\end{theorem} 
\noindent \textbf{Proof}   	From (\ref{gxr}), the GCRJ can be written as
     \begin{align}\label{cc}
     G \xi_m J(X) &=-\frac{1}{2} \int_{0}^{1} \frac{u^m}{f_X(F_X^{-1} (1-u))}du, 
     \end{align}
   and from (\ref{gxp}) the GCPJ can be written as
     	\begin{align}\label{dd}
      	G\bar{\xi}_m J(X) &=-\frac{1}{2} \int_{0}^{1} \frac{u^m}{f_X(F_X^{-1} (u))}du.
     \end{align}
     When $G\bar {\xi}_m (J(X))=G\xi_m(J(X))$ holds then from (\ref{cc}) and (\ref{dd}), 
we have\[\int_{0}^{1} \eta(u) u^m du =0,\] where 
     \begin{align}
     	 \eta(u)=\frac{1}{f_X(F_X^{-1} (1-u))}-\frac{1}{f_X(F_X^{-1} (u))}.
     \end{align}
   Then since $\eta(1-u)=-\eta(u)$, therefore we can write 
   \begin{align} \label{ee}
   	  \int_{0}^{\frac{1}{2}} \eta(u) \left(u^m-(1-u)^m\right)du=0 .
   \end{align}
  Clearly $u^m\leq(1-u)^m $ for all $0\leq u \leq \frac{1}{2}$. Since by assumption $F_X\in \mathbb{C}$, hence from (\ref{ee}), we have $f_X(F_X^{-1} (u))=f_X(F_X^{-1} (1-u))$, so lemma \ref{fa2012} completes the proof. $\blacksquare$

     Generalized cumulative past extropy coresponding to $n$th  upper $k$-record value $	G\bar {\xi}_m (J(U_{n,k}))$ and  Generalized Cumulative Residual extropy coresponding to $n$th $k$ lower record value $G\xi_m(J(L_{n,k}))$ are, respectively, given by
    \begin{align}
    	G\bar {\xi}_m (J(L_{n,k}))&=-\frac{1}{2}\int_{0}^{\infty}{F_X}^m_{L_{n,k}} (x)dx, \label{lnkkcdf}\\
    \text{and} \ \ 	G\xi_m(J(U_{n,k}))&=-\frac{1}{2}\int_{0}^{\infty}\bar{F}^m_{U_{n,k}} (x)dx.\label{unkkcdf}
    \end{align}

The following theorem gives relationship of symmetric distribution with GCPJ and GCRJ of k-record value.
    \begin{theorem}  Let $X_1,X_2,\ldots$ be a random sample of continuous random variables from a population $X$ having cdf $F_X$ and pdf $f_X$.
Then following two statements are equivalent for any $F_X\in \mathbb{C}$:
\begin{enumerate}[(i)]   	
    \item   random variable $X$ has a symmetric distribution;
 \item for a fixed $k\geq 1$,   $G\bar {\xi}_m (J(L_{n,k})) = G\xi_m(J(U_{n,k}))$ for all $n\geq 1$ and for a given fixed number $m \geq 1.$
\end{enumerate}
      \end{theorem}
      \noindent \textbf{Proof} \
      From (\ref{a}) and (\ref{unkkcdf}), the GCRJ of $U_{n,k}$ is 
      \begin{align}\label{gcrex}
      G\xi_m J(U_{n,k})=-\frac{1}{2} \int_{0}^{1} \left( \sum_{i=0}^{n-1} \frac{(-k \log u)^i}{i!} \right) ^m \frac{u^{mk}}{f_X(F_X^{-1}(1-u)))}du  .    
    \end{align}
From (\ref{b}) and (\ref{lnkkcdf}), the GCPJ of $L_{n,k}$ is 
\begin{align}\label{gcpex}
	G\bar{\xi}_m J(L_{n,k})	&=-\frac{1}{2} \int_{0}^{1} \left( \sum_{i=0}^{n-1} \frac{(-k \log u)^i}{i!} \right) ^m \frac{u^{mk}}{f_X(F_X^{-1}(u)))}du.
      \end{align}
     When G$\bar{\xi}_m J(L_{n,k})=G \xi_m J(U_{n,k})$ then from (\ref{gcrex}) and (\ref{gcpex}) we get,
      \begin{align}\label{zz}
      	 \int_{0}^{1} \phi_n^m(u) \eta(u)du=0,       
      \end{align}
    
     \begin{align*}
      \text{ where }\ \ \	\eta(u)=\frac{1}{f_X(F_X^{-1} (1-u))}-\frac{1}{f_X(F_X^{-1} (u))},\\ \text{and} \ \ 
      	\phi_n(u)= u^k\sum_{i=0}^{n-1} \frac{(-k \log u)^i}{i!}.
      \end{align*}
     Further simplification of (\ref{zz}) gives
    \begin{equation}\label{glnx} \int_{0}^{\frac{1}{2}} [\phi_n^m(u)-\phi_n^m(1-u)] \eta(u)du =0. \end{equation}
Note that it is proved in Theroem \ref{thm2}, $\phi_n(u)$ is positive and increasing in $u$. 
Since by assumption $F_X\in \mathbb{C}$ and $\phi_n(u)$ is positive and increasing in $u$, hence from (\ref{glnx}) we have $f_X(F_X^{-1} (u))=f_X(F_X^{-1} (1-u))$, so lemma \ref{fa2012} completes the proof.$\blacksquare$

 \section{Results based on measure of inaccuracy }

     Kerridge's inaccuracy was expanded to the case of a continuous situation by Nath (1968), who also discussed some of their properties. Let $X$ and $Y$ be two continuous random variables with cdf’s $F$ and $G$, and pdf's $f_X(x)$ and $g_Y(x)$, respectively. If $F_X(x)$ is the actual distribution function corresponding to the observations and $G_Y(x)$ is the distribution assigned by the experimenter then the inaccuracy extropy measure (KIJ) is defined as
	\begin{align}\label{31}
		KIJ(X,Y)=-\frac{1}{2} \int f_X(x)g_Y(x)dx.
	\end{align}
   We point out that when $g_Y(x)=f_X(x)$ for all $x$, (\ref{31}) becomes extropy. The $KIJ(U_{n,k},X)$ and $KIJ(L_{n,k},X)$ are given as 
    \begin{align*}
    	KIJ(U_{n,k},X)&=-\frac{1}{2} \int_{S_X} f_{U_{n,k}} (x)f_X(x)dx \\
    	             &= -\frac{1}{2} \int_{0}^{1} \frac{(-klog(u))^{n-1}}{(n-1)!} ku^{k-1} f_X(F_X^{-1} (1-u)) du,  \\ \text{and} \ \
    	KIJ(L_{n,k},X)&=-\frac{1}{2} \int_{S_X} f_{L_{n,k}} (x)f_X(x)dx\\
    	             &= -\frac{1}{2} \int_{0}^{1} \frac{(-klog(u))^{n-1}}{(n-1)!} ku^{k-1} f_X(F_X^{-1} (u)) du, \\
    	             \end{align*}
 respectively.
 The following theorem establish relationship of symmetric distribution with $ KIJ(U_{n,1},X)$ and $ KIJ(L_{n,1},X) $.
   \begin{theorem}   The following two statements are equivalent for any $F_X\in \mathbb{C}$:
 \begin{enumerate}[(i)] \item random variable $X$ has a symmetric distribution;
\item $ KIJ(U_{n,1},X)=KIJ(L_{n,1},X) $. \end{enumerate}  
	\end{theorem} 

\noindent \textbf{Proof} If $ KIJ(U_{n,1},X)=KIJ(L_{n,1},X) $ then 
\begin{align}\label{p31}
	\int_{0}^{1} (-log(u))^{n-1} [f_X(F_X^{-1}(u))-f_X(F_X^{-1}(1-u))]du=0 .
\end{align}
Further simplification gives
\begin{align}\label{p31}
	\int_{0}^{\frac{1}{2}} \left[(-log(u))^{n-1}-(-log(1-u))^{n-1} \right] [f_X(F_X^{-1}(u))-f_X(F_X^{-1}(1-u))]du=0.
\end{align}
Since by assumption $F_X \in \mathbb{C}$ and $(-log(u))^{n-1}$ is positive and strictly decreasing function\ for $u \in (0,\frac{1}{2})$. This implies 
	$f_X(F_X^{-1}(u))=f_X(F_X^{-1}(1-u)) $	for almost all $u  \in  (0,\frac{1}{2}).$ Hence the result by lemma \refeq{fa2012}.$\blacksquare$

 The cumulative residual inaccuracy (CRI) and cumulative past inaccuracy (CPI) measures, which are extensions of the corresponding cumulative residual entropy and cumulative past entropy, respectively, were studied by Kundu et al.(2016). Hashempour and Mohammadi (2022) introduced a new measure of uncertainty that will be called cumulative past extropy inaccuracy measure. Cumulative residual extropy inaccuracy (CRIJ) and cumulative past extropy inaccuracy (CPIJ) measures, which are extensions of the corresponding cumulative residual extropy and Cumulative past extropy, respectively, are defined as :
 \begin{align}
 	CRIJ(X,Y)&=-\frac{1}{2} \int \bar{F}_{X} (x)\bar{F}_Y(x)dx, \\
 \text{and} \ \	
 CPIJ(X,Y)&=-\frac{1}{2} \int F_{X} (x)F_Y(x)dx.
 \end{align}
The CRIJ of $U_{n,k}$ and $X$ is given as, 
\begin{align}
	CRIJ(U_{n,k}, X)&=-\frac{1}{2} \int \bar{F}_{U_{n,k}} (x)\bar{F}_X(x)dx \nonumber\\
	&=-\frac{1}{2} \int_{0}^{\infty} \bar{F}^k_X (x)   \sum_{i=0}^{n-1} \frac{(-k \log\bar{F}_X (x) )^i}{i!} \bar{F}_X(x)dx \nonumber\\
	&=-\frac{1}{2} \int_{0}^{1} u^{k+1} \sum_{i=0}^{n-1} \frac{(-k \log u) )^i}{i!} \frac{du}{f_X({F_X}^{-1}(1-u))}. \label{cri}
\end{align}
 The CPIJ of $L_{n,k}$ and $X$ is given as,
 \begin{align}
 	CPIJ(L_{n,k},X)&=-\frac{1}{2} \int F_{L_{n,k}} (x)F_X(x)dx \nonumber\\
 	              &=-\frac{1}{2} \int_{0}^{\infty} {F_X}^k_X (x)   \sum_{i=0}^{n-1} \frac{(-k \log F_X (x) )^i}{i!} F_X(x)dx \nonumber \\
 	              &=-\frac{1}{2} \int_{0}^{1} u^{k+1} \sum_{i=0}^{n-1} \frac{(-k \log u) )^i}{i!} \frac{du}{f_X(F_X^{-1}(u))}. \label{cpi}
 	            \end{align}
The following theorem establish relationship of symmetric distribution along with $ CRIJ(U_{n,k},X) $ and $ CPIJ(L_{n,k},X)$.
\begin{theorem}  Let $X_1,X_2,\ldots$ be a random sample of continuous random variables from a population $X$ having cdf $F_X$ and pdf $f_X$.
Then following two statements are equivalent for any $F_X\in \mathbb{C}$:
\begin{enumerate}[(i)]   	
    \item   random variable $X$ has a symmetric distribution;
 \item for a fixed $k\geq 1$,   $ CRIJ(U_{n,k},X)=CPIJ(L_{n,k},X) $ for all $n\geq 1$.
\end{enumerate}
      \end{theorem}
 \noindent \textbf{Proof} When  $ CRIJ(U_{n,k},X)=CPIJ(L_{n,k},X) $ then   after simlifications, we get
 \begin{align*}
 	\int_{0}^{1} \phi(u) \eta(u)du =0,
 \end{align*}	
 where $ \phi(u) $ and $\eta(u)$, respectively, are  \\
 \begin{align*}
 	&\phi(u)= u^{k+1} \left( \sum_{i=0}^{n-1} \frac{(-k \log u)^i}{i!} \right),  \\ \text{and} \ \  	&\eta(u)=\frac{1}{f_X(F_X^{-1} (1-u))}-\frac{1}{f_X(F_X^{-1} (u))}.
 \end{align*}
 Since $\eta(1-u)=-\eta(u)$, integral reduces to 
 \begin{align}\label{llnx}
 	\int_{0}^{\frac{1}{2}} (\phi(u)-\phi(1-u) )\eta(u)du =0.
 \end{align}
 Note that $\phi(u)\geq 0$. Also
\begin{align*}
\frac{d\phi(u)}{du}&=ku^{k}\left(\sum_{i=0}^{n-1} \frac{(-k \log u)^i}{i!}-\sum_{i=0}^{n-2} \frac{(-k \log u)^i}{i!}\right)+u^{k}\sum_{i=0}^{n-1} \frac{(-k \log u)^i}{i!}\\
&\geq 0.
\end{align*}    
Since by assumption $F_X\in \mathbb{C}$ and $\phi(u)$ is positive and increasing in $u$, hence from (\ref{llnx}), we have $f_X(F_X^{-1} (u))=f_X(F_X^{-1} (1-u))$, so the lemma \ref{fa2012} completes the proof. $\blacksquare$\\
 
\section{Examples}
To illustrate the results from earlier sections, we now look at a few examples. For similar entropy-based characteristics and examples  see Ahmadi (2021).
\subsection{ Power function distribution :} The pdf and cdf respectively of power function distribution are
\begin{align}
	f_X(x)=\theta x^{\theta -1} \ \ \ \text{and} \ \ \ F_X(x)= x^\theta, \ \  0<x<1, \ \theta >0. 
\end{align}
Using above pdf and cdf, 
\begin{align*}
	f_X(F_X^{-1}(u))=\theta u^{\frac{\theta -1}{\theta}} \ \ \text{and} \ \ f_X(F_X^{-1}(1-u))=\theta (1-u)^{\frac{\theta -1}{\theta}}, \ \ 0<u<1.
\end{align*}
  
  If $\theta >1$ then $f_X(F_X^{-1}(u)) \leq f_X(F_X^{-1}(1-u))$ for u $\in$ (0,$\frac{1}{2}$ ); if $\theta <1$ then $f_X(F_X^{-1}(u)) \geq f_X(F_X^{-1}(1-u))$ for u $\in$ (0,$\frac{1}{2}$) and if $\theta=1$ then $f_X(F_X^{-1}(u)) = f_X(F_X^{-1}(1-u))$ for u $\in$ (0,$\frac{1}{2}$).  Thus power distribution belongs to class $\mathbb{C}$.\\
  \indent Consider Theorem \ref{g} for a moment. The difference
  \begin{align*}
  	\Delta_1 &= \xi J(X) -\bar{\xi} J(X) \\
  	 &=-\frac{1}{2} \int_{0}^{\frac{1}{2}} \left( \frac{1}{f_X(F_X^{-1} (1-u))}-\frac{1}{f_X(F_X^{-1} (u))}\right) (2u-1)du\\
  	 &=-\frac{1}{2\theta} \int_{0}^{\frac{1}{2}} \left( \frac{1}{(1-u)^{\frac{\theta-1}{\theta}}}-\frac{1}{u^{\frac{\theta-1}{\theta}}}\right) (2u-1)du\\
  	 &=\frac{1-\theta}{2(\theta +1)}.
    \end{align*}
   Then for $u \in (0,\frac{1}{2}$),\ \ $\Delta_1=0$  for $ \theta =1,$ hence the power distribution is symmetric for $\theta=1$.\\
   \indent   One can examine other theorems in similar fashion for Power function distribution.

   \subsection{ Pareto distribution :} The pdf and cdf respectively of Pareto distribution are
    \begin{align}
   	f_X(x)=\theta x^{-\theta -1} \ \ \ \text{and} \ \ \ F_X(x)= 1-x^{-\theta}, \ \  x>1 , \ \theta >0. 
   \end{align}
 Using above pdf and cdf, 
\begin{align*}
	f_X(F_X^{-1}(u))=\theta (1-u)^{\frac{\theta +1}{\theta}} \ \ \text{and} \ \ f_X(F_X^{-1}(1-u))=\theta (u)^{\frac{\theta +1}{\theta}}, \ \ 0<u<1.
\end{align*}
   For $ u\in (0,\frac{1}{2}),\ f_X(F_X^{-1}(u)) > f_X(F_X^{-1}(1-u))$, that is, Pareto distribution belongs to class $\mathbb{C}$.\\
   \indent Consider Theorem \ref{thm2} for a moment. The difference
   	\begin{align*}
   		\Delta_2 &= \xi J(U_{n,k}) -\bar{\xi} J(L_{n,k}) \\
   		&=-\frac{1}{2} \int_{0}^{\frac{1}{2}} \left( \phi^2(u)-\phi^2(1-u) \right) \left( \frac{1}{f_X(F_X^{-1} (1-u))}-\frac{1}{f_X(F_X^{-1} (u))}\right)du\\
   		&=-\frac{1}{2\theta} \int_{0}^{\frac{1}{2}} \left( \phi^2(u)-\phi^2(1-u) \right) \left( \frac{1}{ (u)^{\frac{\theta +1}{\theta}}}-\frac{1}{ (1-u)^{\frac{\theta +1}{\theta}}}\right)du
   	\end{align*}
   where $\phi_n(u)= u^k\sum_{i=0}^{n-1} \frac{(-k \log u)^i}{i!}$ is positive and increasing function of $u$ as shown in Theorem \ref{thm2}. 
   	Then for $u \in (0,\frac{1}{2}$), it is easy to show that $\Delta_2>0$ for all $\theta >0.$ This is due to the fact that the Pareto distribution is not symmetric.\\
   	\indent One can examine other theorems in similar fashion for Pareto distribution.
   	\subsection{ Exponential distribution :}
   	The pdf and cdf respectively of exponential distribution with mean 1 are
   	\begin{align}
   		f_X(x)=e^{-x} \ \ \ \text{and} \ \ \ F_X(x)= 1-e^{-x}, \ \  x>0. 
   	\end{align}
   	Using above pdf and cdf, 
   	\begin{align*}
   		f_X(F_X^{-1}(u))=1-u \ \ \text{and} \ \ f_X(F_X^{-1}(1-u))=u, \ \ 0<u<1.
   	\end{align*}
   	For $ u\in (0,\frac{1}{2}),$ $$f_X(F_X^{-1}(u)) > f_X(F_X^{-1}(1-u)),$$ \text{and} for $ u\in (\frac{1}{2},1),$\ $$f_X(F_X^{-1}(u)) < f_X(F_X^{-1}(1-u)),$$ that is, exponential distribution belongs to class $\mathbb{C}$. Consider Theorem \ref{thm3} for a moment. The difference
   	\begin{align*}
   		\Delta_3 &= G\bar {\xi}_m (J(X))-G\xi_m(J(X))\\
   		&=-\frac{1}{2} \int_{0}^{\frac{1}{2}} \left[ \frac{1}{f_X(F_X^{-1} (1-u))}-\frac{1}{f_X(F_X^{-1} (u))} \right] \left(u^m-(1-u)^m\right)du\\
   		&=-\frac{1}{2}  \int_{0}^{\frac{1}{2}} \frac{(1-2u)}{u(1-u)} \left(u^m-(1-u)^m\right)du
   	\end{align*}
   It is easy to show that $\Delta_3 >0$ for all $u\in (0,\frac{1}{2})$ and for all $m \geq 1.$ This is due to the fact that exponetial distribution is not symmetric.\\
  \indent One can examine other theorems in similar fashion for exponential distribution.
  \subsection{ Uniform distribution :} 	The pdf and cdf, respectively of uniform distribution on $(0,1)$ are
  \begin{align}
  	f_X(x)=1 \ \ \ \text{and} \ \ \ F_X(x)= x, \ \  x\in (0,1).
  \end{align}
 
  \indent For $ u\in (0,1), \ F_X^{-1}(u)=u$. Then $f_X(F_X^{-1}(u))=  f_X(F_X^{-1}(1-u)) $. As a result, the equality is true in $\mathbb{C}$ since the uniform distribution is symmetric and all of the established theorems in the previous section are true for the uniform distribution.            \\
   	\subsection{ Standard Normal distribution :}
   	The pdf and cdf, respectively of standard Normal distribution are
   	\begin{align}
   		f_X(x)=\frac{1}{\sqrt{2\pi}}e^{-\frac{x^2}{2}} \ \ \ \text{and} \ \ \ F_X(x)= \frac{1}{2} \left[ 1+erf(\frac{x}{\sqrt{2}})\right], \ \  x\in \mathbb{R},
   	\end{align}
   	 where erf(.) is error function given by $$erf(x)=\frac{2}{\sqrt{\pi}} \int_{0}^{x}e^{-t^2}dt$$ which is an odd function. Also note that $f_X(x)$ is even function. \\
   	 \indent For $ u\in (0,1),\  F_X^{-1}(u)=\sqrt{2}\ erf^{-1}(2u-1)$. Then $f_X(F_X^{-1}(u))=  f_X(F_X^{-1}(1-u)) $.  As a result, the equality is true in $\mathbb{C}$ since the standard normal distribution is symmetric and all of the established theorems in the previous section are true for the standard normal distribution.    \\
   	 
   	 \section{Conclusion and future work}
   	 Our purpose of this study was to investigate and provide some new characteristics and results of continuous symmetric distribution based on extropy. For Similar characteristics results for continuous symmetric distribution based on entropy, see Ahmadi (2021).\\
   	 \indent Xiong et al. (2021) used one characterization of symmetric distribution based on extropy of $k$th upper record value and $k$th lower record value in testing symmetry of probability distribution. Also, Jose and Sather (2022) used one characterization of symmetric distribution based on extropy of $n$th upper $k$-record value and $n$th lower $k$-record value in testing the symmetry of probability distribution. Therefore, one can use the characterization result provided in this paper to construct new test statistics in testing the symmetry of continuous probability distribution.\\
   	 \indent For example, we may be interested in testing $H_0: f_X(k+x)=f_X(k-x)$  for all $x\in S_X $ against $H_1: f_X(k+x) \neq f_X(k-x) $ for some $x\in S_X $.  Based on theorem \ref{g}, equivalent test is $H_0: \bar{\xi} J(X)-\xi J(X)=0$   against $H_1: \bar{\xi} J(X)-\xi J(X) \neq 0$. Thus, one can get test statistics as an estimate of  $\Delta = \bar{\xi} J(X)-\xi J(X) $ in testing symmetry. Empirical estimators of $\bar{\xi} J(X)$ and $\xi J(X)$ can be used that have been studied by researchers, see Jose and Sathar (2022), Xiong et al.(2021), Park (2020), Zardasht et al (2015), Noughabi (2015). Thorough research is needed to work on this problem.
   	     \\
   	     \\
   	     \\
    \textbf{ \Large Funding} \\
        \\
	Santosh Kumar Chaudhary would like to thank Council Of Scientific And Industrial Research (CSIR), Government of India ( File Number 09/0081(14002)/2022-
	EMR-I ) for financial assistance.\\
	\\
	\textbf{ \Large Conflict of interest} \\
	\\
	The authors declare no conflict of interest.\\


\vspace{.5in}
\noindent
{\bf  Nitin Gupta} \\
Department of Mathematics,\\
Indian Institute of Technology Kharagpur\\
Kharagpur-721302, INDIA\\
E-mail: nitin.gupta@maths.iitkgp.ac.in\\

\vspace{.1in}

\noindent
{\bf Santosh Kumar Chaudhary}\\
Department of Mathematics,\\
Indian Institute of Technology Kharagpur\\
Kharagpur-721302, INDIA\\
E-mail: skchaudhary1994@kgpian.iitkgp.ac.in \\

	\end{document}